\documentclass[12pt,letterpaper]{iopart}
\usepackage{iopams}
\newcommand{\ignore}[1]{}
\def\calm{{\cal M}}
\def\caln{{\cal N}}
\def\udel{\Delta}
\def\del{\delta}
\def\ep{\epsilon}
\def\what{\widehat}
\def\wtilde{\widetilde}
\def\tilden{\widetilde{\caln}}
\def\ugam{\Gamma}
\def\pf{\noindent{\bf Proof.}~}
\newtheorem{lem}{Lemma}
\newtheorem{thm}{Theorem}
\newtheorem{rem}{Remark}
\newtheorem{dfn}{Definition}
\begin{document}
\title{Escape from a circle and Riemann hypotheses}
\author{Leonid A Bunimovich$^1$ and Carl P Dettmann$^2$}
\address{$^1$ School of Mathematics, Georgia Institute of 
Technology, Atlanta, GA 30332, USA}
\address{$^2$ School of Mathematics, University of Bristol,
Bristol BS81TW, United Kingdom}

\ead{bunimovh@math.gatech.edu}
\begin{abstract}We consider open circular billiards with
one and with two holes.  The exact formulas for escape
are obtained which involve the Riemann zeta function and 
Dirichlet $L$ functions. It is shown that the problem of finding the
exact asymptotics in the small hole limit for escape
in some of these billiards is equivalent to the Riemann hypothesis.
\end{abstract}

\maketitle

\baselineskip=22pt
\section{Introduction}

The theory of open dynamical systems where orbits may disappear upon
reaching some region in phase space just recently started to attract the
attention of mathematicians.  A natural reason for
this is because it is much harder in general to study open systems in
comparison with closed dynamical systems.  The theory of closed
dynamical systems is rather well developed.  Therefore these
relatively few mathematical studies of open dynamical systems
are naturally based on ideas and techniques developed in dealing with 
closed systems.  However, the studies of open systems may also bring
new insights to the understanding of dynamics of closed systems as well.

One of possible paths to explore in this direction was recently
suggested in Ref.~\cite{1}.  Let us consider open systems with several ``holes"
(regions where orbits disappear upon hitting them) and compare behaviors
of such ``many-holes" systems with the corresponding single-hole
systems.  The idea is that such comparison may shed a light on
understanding dynamics of a closed system one gets by 
``patching" all the holes
in the open systems.  (This approach could be of interest for
geometric theory of dynamical systems by studying e.g., geodesic flows
on manifolds with ``holes.")

This idea has also a lot of potential applications for the real world
``physical'' systems~\cite{1}.  Indeed, in experimental studies 
researchers perform measurements outside a region of interest (``container") by
e.g., measuring fluxes out of the container~\cite{2,3,4}.  
Actually our approach arose from the claim made by one of the authors, who
was inspired by experiments with the optical billiards, that a comparison of
escape rates through one and through two holes may shed some light 
on the dynamics of the corresponding closed systems.
This approach
may even have some potential industrial applications related, e.g., to
the optimal placement of the holes in order to maximize (or minimize)
corresponding fluxes.

One of the most natural classes of the open systems to study is formed
by the open billiards.  The first application of the one-many
holes interplay idea has already led to a very surprising and,
in a sense, remarkable result.  Seemingly the simplest problem of 
this type, namely the comparison of dynamics of the circular billiards
with one and with two opposite holes, is equivalent to the Riemann hypothesis
(RH) \cite{1}.  Formally, this result opens up a possibility to verify RH in
real physical experiments.  (It does not seem however that it is going to be
very fruitful and, especially, practical approach because the modern
computers provide more efficient tools for that in numerical experiments.)
On the other hand, this result demonstrates that indeed the studies
of open systems (and particularly of open billiards) may bring about some
interesting and unexpected advances in traditional areas of research.

The purpose of this paper is to provide the proofs of the results
announced in \cite{1}.

\section{Definitions and notations}
Consider a billiard on the unit disk $D$ i.e., a dynamical system generated
by the motion of a point particle with a constant speed within $D$ with
elastic collisions (angle of incidence equals angle of reflection)
from its boundary.  Without any loss of generality we assume that the 
particle's speed is identically one, and therefore its velocity is
completely defined by an angle 
$\vartheta$ it makes with the horizontal direction,
$-\pi<\vartheta\le\pi$.  Billiards are Hamiltonian systems.  Therefore
the Liouvillean measure (the phase volume) in the phase space $\calm$ is
preserved under the dynamics (a billiard flow) $\{S^t\}$, $-\infty<t<\infty$.

Let $M=\{(\beta,\psi):-\pi<\beta\le\pi, -\frac\pi2\le\psi\le
\frac\pi2\}$, where $\beta\in \partial D$.  The billiard flow in
$D$ induces the billiard map $T:M\to M$ defined as 
\begin{equation}\label{eq1}
T(\beta,\psi)=
(\beta +\pi-2\psi,\psi)\end{equation} 
where $\psi$ is the angle between the outward trajectory and the inner normal
at $\beta\in\partial D$, and all the angles in \eref{eq1} are taken
modulo $2\pi$.  The natural projection $M\to\partial D$ we denote
by proj.  Thus proj$(\beta,\psi)=\beta$, where $(\beta,\psi)\in M$.

It is well known that orbits of the billiard in a circle are either 
periodic with a period $q$ (if $\psi=\frac\pi2-\frac pq\, \pi$, where
$p$ and $q$ are co-prime integers and $p<q$) or everywhere dense
in $\partial D$ (if $\psi$ is incommensurable with $\pi$).  The billiard
map $T$ preserves the measure 
\begin{equation}\label{eq2}
d\mu=\frac{\cos\psi}{4\pi}\, d\psi\, d\beta\end{equation}
which is the projection of the Liouvillean measure
in the phase space $\calm$ of the billiard flow onto $M$.

Suppose that two (possibly overlapping) holes
$H_1=\{\beta,0<\beta<\udel\}$ and $H_2=\{\beta:
\theta<\beta<\theta+\udel\}$, $0\le\theta<2\pi$, $\udel >0$,
are placed at the boundary $\partial D$.
Consider now a new dynamical system, an open billiard in $\partial D$
with holes $H_1$ and $H_2$.  In this open billiard any orbit
$(\beta_0,\psi_0)$ moves under the billiard map \eref{eq1} until
it hits one of the holes $H_1$  and $H_2$.  When the orbit hits
$H_1\cup H_2$ it ``disappears" (escapes).

Obviously, almost all (with respect to the measure $\mu$) orbits
will eventually escape.  The only orbits that never escape are such
periodic orbits that never hit $H_1\cup H_2$.  Denote
$\what H_i=\{(\beta,\psi):\beta\in H_i\}$, $i=1,2$.  Thus 
proj~$\what H_i=H_i$, $i=1,2$. By dist$(\beta_1,\beta_2)$,
$\beta_i\in\partial D $, $i=1,2$, we will mean the length
of the shortest arc between $\beta_1$ and $\beta_2$.

Let $N(\beta_0,\psi_0)$, $(\beta_0,\psi_0)\in M$ be a (minimal) number
of reflections from the boundary after which the orbit 
$T^n(\beta_0,\psi_0) = (\beta_n,\psi_n)$, $n=1,2,\dots$ escapes from
the circle.  (If the orbit of $(\beta_0,\psi_0)$ never escapes we
set $N(\beta_0,\psi_0)=\infty$.)  Therefore the orbit $(\beta_0,\psi_0)$
escapes from the circle in a (real, continuous) time 
$\tau(\beta_0,\psi_0)=2\cos\psi_0 N(\beta_0,\psi_0)$.  

\section{Structure of the set of orbits not escaping in time $t$}
Clearly, the only orbits that never escape are those periodic orbits
with periods $q<\frac{2\pi}\udel$ which never hit the holes
$H_1\cup H_2$.

All orbits of a rotation of the circle at any irrational (with respect
to $\pi$) angle are everywhere dense (on the circle).  Therefore all 
orbits of irrational rotations of the circle eventually escape 
regardless of the size $\udel >0$ of the holes.

We will always assume in what follows that $\udel < \frac\pi2$.  Then there
exists a periodic orbit of period two which never escapes.  

The following statement obviously holds.

\begin{lem}\label{lem1}
If $\udel < \frac\pi n$, $n\ge 2$, then there exists a periodic orbit
of period $n$ which never escapes.
\end{lem}

\begin{lem}\label{lem2new} 
If $\udel < \frac\pi n$, $n\ge 2$, then for any $t>0$ there exists a 
nonperiodic orbit which does not escape till time $t$.
\end{lem}

\pf Let $\frac\pi n-\udel =\del>0$.  By Lemma~\ref{lem1} there exists a 
periodic orbit $(\beta_i,\hat \psi)$, $i=1,2,\dots, n$, $T(\beta_i,\hat\psi)
=(\beta_{i+1},\hat\psi)$ if $1\le i\le n-1$, $T(\beta_n,\hat\psi)=
(\beta_1,\hat\psi)$, such that $\min_{1\le j\le n}
\mbox{\rm dist} (\beta_j,\{H_1\cup H_2\})=\frac\del 2$.  
Therefore for any $t>0$ and for any $j$, $1\le j\le n$, there exists
such $\alpha_{j,t}=\alpha_{j,t}(\del)$ that orbits of all points
$(\beta_j,\psi)$ with $|\psi-\hat\psi|<\alpha_{j,t}$ do not escape
from the circle till time $t$.

Clearly, smaller the difference $|\psi-\hat\psi|$ is, longer the 
corresponding orbit will not escape.

Denote by $\caln_t$ the set of all orbits that do not escape till
the time $t$.  We will show that for sufficiently large $t$, the set 
$\caln_t$ can be decomposed into the union of nonintersecting neighborhoods
of never escaping periodic orbits.

\begin{lem}\label{lem2} 
Let $x'=(\beta, \psi')$ and $x''=(\beta,\psi'')$ be the points 
of two never escaping periodic orbits with periods $n'$ and $n''$ 
respectively.  Then $x'$ and $x''$ belong to different connected
components of the sets $\caln_t$ if $t>2(n'+n''-1)\max(\cos\psi',\cos\psi'')$.
\end{lem}

\pf Suppose $n''>n'$ and $n''\ne kn'$, where $k>0$ is an
integer.
Consider the set $A_{\beta,\psi',\psi''}=
\{(\beta,\psi):\psi'\le \psi\le \psi''\}$.  (We assumed here that 
$\psi''>\psi'$.  But the same argument can be applied if
$\psi''<\psi'$.)  Clearly proj$(T A_{\beta,\psi',\psi''})$ is the
arc of $\partial D$ between the points proj$(Tx')$ and 
proj$(Tx'')$.  Analogously proj$(T^mA_{\beta,\psi',\psi''})$,
$1\le m\le n'$, is the arc between proj$(T^mx')$ and 
proj$(T^mx'')$.  

Now we show that  
proj$\left(\bigcup^{n'+n''-1}_{m=1} T^mA_{\beta,\psi',\psi''}\right)=
\partial D$.  Recall that $x'$ and $x''$ are never escaping periodic orbits.
Observe that the circle $\partial D$ is divided by the points
proj$(T^mx')$, $0\le m\le n'-1$, and proj$(T^mx'')$, $0\le m\le n''-1$ 
into at most $n'+n''-1$ arcs.  We will call these arcs minimal arcs.  
Clearly each arc between proj$(T^mx')$ and proj$(T^mx'')$, covers
at least one minimal arc that was not covered by the arcs between
proj$(T^{m-1}x')$ and proj$(T^{m-1}x'')$.  Indeed, otherwise the
corresponding preimage of this arc also covered only minimal arcs
already covered by the previous iterations of $T$.  Continuing
this argument we will come back to the point $\beta\in \partial D$
and to a contradiction.  Therefore, not more than by $(n'+n''-1)$
iterations all the circle $\partial D$ will be covered.  
Hence, within $(n'+n''-1)$ iterates of the billiard map
$T$ the projections on $\partial D$ of the corresponding images of the set
$A_{\beta,\psi',\psi''}$ will completely cover the hole $H_1$ and
the hole $H_2$.  Thus after that time the iterates of $x'$ and 
$x''$ will belong to different connected components of the sets 
$\caln_t$.
The case $n''=kn'$ is even simpler to consider and
the analysis goes along the same lines.

We will show now that any connected component of the set
$\caln_t$ contains some interval 
$I_{\psi_\frac mn,\beta_1,\beta_2}=\{(\beta,\psi):\psi=\frac\pi 2
-\frac mn\, \pi$, $\beta_1\le\beta\le \beta_2\}$.  Observe that the 
orbit of the point $(\beta,\frac\pi2-\frac mn\, \pi)$ has period $n$
and the $\beta$ values are equally spaced at intervals of
$\frac{2\pi}n$.

The next statement assures that any connected component
of the set of nonescaping orbits $\caln_t$ for any $t$ contains
never escaping periodic orbits.  Let $(a,b)$ denote the
gcf~$(a,b)$.

\begin{lem}\label{lem4}
All nonempty connected components of the set $\caln_t$ for all $t>0$ contain a
point of a never escaping periodic orbit.  

\end{lem}

\pf Let $m_i/n_i,\quad(m_i,n_i)=1,\quad i=1,2$ be consecutive fractions among all $m/n,\quad(m,n)=1,\quad \frac{2\pi}n>\udel$, i.e.~Farey numbers.  Assume at first that for some $\beta$ periodic orbits of both points $x_1=(\beta, \frac{\pi}{2}-\frac{m_1}{n_1}\pi),\quad x_2=(\beta,\frac{\pi}{2}-\frac{m_2}{n_2}\pi)$ never escape.

Consider the sets
\begin{eqnarray*}
A^{(1)}&=A^{(1)}_{\beta,\frac\pi2-\frac{m_1}{n_1}\pi,
\frac\pi2 -\frac{m_1}{n_1}\pi + \frac\pi{n_1(n_1+n_2)}}\quad\mbox{and}\\
A^{(2)}&=A^{(2)}_{\beta,\frac\pi2-\frac{m_2}{n_2}\pi
 - \frac\pi{n_2(n_1+n_2)},\frac\pi2-\frac{m_2}{n_2}\pi}
\end{eqnarray*}
Because $\frac{m_1}{n_1}$ and $\frac{m_2}{n_2}$ are the consecutive
Farey numbers (with denominators not exceeding 
$[\frac{2\pi}\udel]$)
$$A_{\beta,\frac\pi2-\frac{m_1}{n_1}\pi,\frac\pi2-
\frac{m_2}{n_2}\pi} =A^{(1)}\cup A^{(2)},$$
where $[a]$ denotes the integer part of the number $a$.

Consider proj$(T^{n_i}A^{(i)})$, $i=1,2$.  It is easy to calculate that
the length of each of these two arcs of $\partial D$ equals
$2\pi/(n_1+n_2)$.  
The well known property of Farey numbers ensures that
$(n_1+n_2)>[\frac{2\pi}\udel]+1$.  Therefore the distance
between the projection $\beta$ of any point $(\beta,\psi)\in A^{(i)}$,
$i=1,2$, and the projection proj$(T^{n_i}(\beta,\psi))$ does not exceed
$\frac{2\pi}{n_1+n_2}<\udel$.

The last inequality ensures that all images of the point $(\beta,\psi)
\in A^{(1)}\cup A^{(2)}$
belong (before escaping) to the same connected component of the sets
$\caln_t$ to which belong at least one point of the periodic
orbit $(\beta,\frac\pi 2-\frac{m_1}{n_1}\pi)$ or of
$(\beta,\frac\pi2-\frac{m_2}{n_2}\pi)$ or of each of these periodic
orbits (if $t$ is sufficiently small).

Let now assume that one or both of the points $x_i,\quad i=1,2$
do escape. Then for some $k>0$ proj$(T^kx_1)$ belongs to $H_1\cup H_2$.
(The consideration is absolutely analogous if the orbit of $x_2$ escapes
at some bounce off the boundary of the disk $D$.) Consider all the points
of the set $T^kA^{(1)}$ that have not escaped to this moment. Among these
points must be one with projection at the end point of the hole $H_i$,
where $T^kx_1$ escaped. Consider now the point $x_3$ that goes from this 
point under the angle $\frac{\pi}2 - \frac{m_1}{n_1}\pi$. Obviously this point is periodic
and the corresponding periodic orbit will never escape.  

Consider all the points of the set $T^kA^{(1)}$ that have not escaped to 
this moment. Clearly, the length of the projection of this set onto the 
boundary of $D$ does not exceed $\pi/(n_1+n_2)<\Delta/2$. Therefore all periodic 
orbits that go from the points of this set under the angle $\pi/2 - 
\pi m_1/n_1$ will never escape.

On the other hand 
from the above argument in the proof it 
follows that all the points of the images of the set $A^{(1)}$ will belong
(until the escape) to the same connected component of the set $\caln_t$
as the images of the never escaping periodic point $T^{-k}x_3$. 
This completes the proof of Lemma 4.

Lemmas~\ref{lem2} and \ref{lem4} imply that the following statement
holds.

\begin{thm}\label{thm1new}
Let $t>4[\frac{2\pi}\udel]$. Then every connected component $B_i$, 
$i=1,2,\ldots,m\quad m=m(\udel)$  of the set 
$\caln_t$ of orbits never escaping till time $t$ contains a unique segment
$I_i= \{(\beta,\psi),\quad\beta_{i,1}<\beta<\beta_{i,2}\}$ consisting of never escaping
periodic orbits.  
\end{thm}

\section{Probability of not escaping till time $t$}

We now compute the measure $\mu(\caln_t)$ of the set of all orbits
that do not escape till time $t$.  Lemmas~\ref{lem1}--\ref{lem4} imply
that $\mu(\caln_t)>0$ for any $t<\infty$.  

Denote $\psi_{m,n} = \frac\pi2-\frac mn\, \pi$, where $m<n$,
$(m,n)=1$, $n<\left[\frac{2\pi}\Delta\right]$.  Clearly
$\caln_t\subset M\backslash \bigcup^{n-1}_{k=0} (T^{-k}
(\hat H_1\cup\hat H_2))$ if $t>2\left[\frac{2\pi}\Delta\right]$.
In what follows we will always assume that $t>\frac{8\pi}\Delta$.

Let $(\beta,\psi)\in \caln_t$.  Then in view of 
Lemma~\ref{lem4} and Theorem~\ref{thm1new} the coordinate
$\psi$ can be uniquely represented as $\psi=\psi_{m,n} + \eta$,
where $|\eta|<\Delta/2$.  It is easy to see that the orbit of the
point $(\beta,\psi)$ escapes not later than at the time $t$ if there
exists such  integer $k$, $0\le k\le\left[ \frac t{2\sin \left(\frac mn\, \pi
+\eta\right)}\right]$ that
\begin{equation}\label{eq3}
\mbox{proj }T^k(\beta,\psi)\in H_1\cup H_2\end{equation}
Denote by $T_{m,n}$ rotation of the circle $\partial D$ on the
angle $\psi_{m,n}$. Lemma~{\rm\ref{lem4}} and Theorem~{\rm\ref{thm1new}}
ensure that every connected component of the set $\caln_t$ can be
uniquely represented as $\caln^t_{\psi_{m,n,j}}$ where 
$m<n$, $(m,n)=1$, $n<\left[ \frac{2\pi}\Delta\right]$ and 
$0\le j\le 2n-1$ if 
$H_2\cap (H_1\cup T_{m,n}H_1\cup T^2_{m,n}H_1\cup\cdots \cup
T^{n-1}_{m,n}H_1)=\emptyset$ or, otherwise, $0\le j\le n-1$.

Let $\theta'=\theta(\bmod\frac{2\pi}n)$.  If $\theta'<\Delta$
then the set $\partial D\backslash\bigcup^{n-1}_{k=0} T^k_{m,n}
(H_1\cup H_2)$ consists of $n$ arcs of the length $\frac{2\pi}n-
\theta'-\Delta$, otherwise it consists of $2n$ arcs, $n$ of
which are of the length $\frac{2\pi}n-\theta'-\Delta$ and another $n$
arcs are of the length $\theta'-\Delta$.  We will call these arcs
complements to a hole's orbit.

Thus one can write
\begin{equation}\label{eq4}
\caln_t=\bigcup_{\stackrel{\stackrel{\scriptstyle m<n}{\scriptstyle (m,n)=1}}
{\scriptstyle n<\left[\frac{2\pi}\udel\right]}}
\bigcup_j
\caln^t_{\psi_{m,n,j}}
\end{equation}
where $0\le j\le 2n-1$ or $0\le j\le n-1$.
Consider now all connected components $\caln^t_{\psi_{m,n,j}}$,
$0\le j\le 2n-1$ (the case when $0\le j\le n-1$ can be treated
analogously and, in fact, is slightly simpler).
According to Lemma~\ref{lem4} and Theorem~\ref{thm1new} all 
$\caln^t_{\psi_{m,n,j}}$ are closed sets for any $m,n,j$.

We will call a connected component $\caln^t_{\psi_{m,n,j}}$ a basic
component if proj$(\caln^t_{\psi_{m,n,j}})$ is adjacent to a hole $H_1$
or $H_2$. (Observe that in the case when $0\le j\le 2n-1$ there are
four basic components for a fixed $m,n$ while if $0\le j\le n-1$
then there are two basic components.)
In each basic component we will single out a closed subset which
will be called a basic set in what follows (instead of its longer
name a basic subset of a basic component).

\begin{dfn}\label{dfn1}
\rm A basic subset $\what\caln^t_{\psi_{m,n,j}}\subset \caln^t_{\psi_{m,n,j}}$
consists of all points $(\beta,\psi)\in\caln^t_{\psi_{m,n,j}}$ such that
$${\rm dist}(\beta,H)\le{\rm dist} \left({\rm proj}\left(
T^n(\beta,\psi)\right),H\right)$$
where $H$ is a hole adjacent to proj$(\caln^t_{\psi_{m,n,j}})$, i.e.,
$H_1$ or $H_2$.
\end{dfn}

Observe that the above inequality is well defined because
dist$(\beta,{\rm proj}(T^n(\beta,\psi))<\Delta$. 
It is easy to see that each basic component $\caln^t_{\psi_{m,n,j}}$ contains
one and only one basic subset.  Therefore for any pair of positive
integers $m<n$, $(m,n)=1$, $\left[\frac{2\pi}n\right]>\Delta$ there
exist four (in the case we consider, otherwise two) basic subsets 
$\what\caln^{t,i}_{m,n}$, $i=1,2,3,4$.

A crucial fact is that for any point $(\beta,\psi)\in \caln^t_{\psi_{m,n}}$
there exist one and only one $k$, $0\le k <n$, and one and only one $i$,
$i=1,2,3,4$, such that $T^k(\beta,\psi)\in\what\caln^{t,i}_{m,n}$.  
Indeed, it follows from the definition of the basic set,
Lemma~\ref{lem2} and the obvious relation that proj$(T^k(\beta,\psi))
\not\in H_1\cup H_2$ for any $0\le k<n$.  

Certainly proj$(\what\caln^{t,i}_{m,n})\cap\left(\bigcup^{n-1}_{k=0}
T^k_{m,n}(H_1\cup H_2)\right)=\emptyset$ for $i=1,2,3,4$.  Therefore projections
of two $(i=1,2)$ basic sets have the length $\frac{2\pi}n -\theta'-\Delta$,
and of two others $(i=3,4)$ $\theta'-\udel$.  (Recall that we consider
only the case $0\le j\le 2n-1$.)

It follows from the definition of the basic sets that if
$(\beta,\psi)\in\what\caln^{t,1}_{m,n}$ then $\psi\ge\psi_{m,n}$ while if
$(\beta,\psi)\in\what\caln^{t,2}_{m,n}$ then $\psi\le\psi_{m,n}$ (or vice
versa).  The same statement is true for the basic sets 
$\what\caln^{t,3}_{m,n}$ and $\what\caln^{t,4}_{m,n}$, i.e., for
$i=3$ and $i=4$.  

Consider now the sets $\wtilde \caln^t_{m,n,k}, 
\skew2\wtilde{\tilden}^t_{m,n,k}\subset \caln^t_{\psi_{m,n}}$,
$k=0,1,\dots, n-1$, where
$$\wtilde\caln^t_{m,n,k} = \left\{(\beta,\psi)\in\caln^t_{\psi_{m,n}}:
T^k(\beta,\psi)\in \what\caln^{t,1}_{m,n}\cup
\what\caln^{t,2}_{m,n}\right\}$$
and
$$\skew2\wtilde{\tilden}^t_{m,n,k} = \left\{(\beta,\psi)\in\caln^t_{\psi_{m,n}}:
T^k(\beta,\psi)\in \what\caln^{t,3}_{m,n}\cup
\what\caln^{t,4}_{m,n}\right\}.$$
It is easy to see that the sets $\wtilde{\caln}^t_{m,n,k}$ and
$\skew2\wtilde{\tilden}^t_{m,n,k}$ are disjoint.  Moreover,
$\wtilde{\caln}^t_{m,n,k_1}\cap \wtilde{\caln}^t_{m,n,k_2}=
\emptyset$
$(\skew2\wtilde{\tilden}^t_{m,n,k_1}\cap\skew2\wtilde{\tilden}^t_{m,n,k_2}
=\emptyset)$ if $k_1\ne k_2$.  Therefore we can write
\begin{equation}\label{eq5}
\caln_t=\bigcup_{\stackrel{\stackrel{\scriptstyle m<n}{\scriptstyle (m,n)=1}}
{\scriptstyle n<\left[\frac{2\pi}\udel\right]}}
\bigcup^{n-1}_{k=0}
\left(\wtilde\caln^t_{m,n,k}\cup
\skew2\wtilde{\tilden}^t_{m,n,k}\right)
\end{equation}
It follows from \eref{eq5} and \eref{eq2} that
\begin{eqnarray}\label{eq6}
\mu(\caln_t)&=
\sum_{\stackrel{\stackrel{\scriptstyle m<n}{\scriptstyle (m,n)=1}}
{\scriptstyle n<\left[\frac{2\pi}\udel\right]}}
\sum^{n-1}_{k=0}
\Biggl(\int^{\frac{2\pi}n-\theta'-\udel}_0
d\beta\int^{\eta^+_{m,n}(\beta,k)}_{-\eta^-_{m,n}(\beta,k)}
\sin(\psi_{m,n}+\eta)d\eta\nonumber\\
&\qquad +
\int^{\theta'-\udel}_0d\beta
\int^{\eta^+_{m,n}(\beta,k)}_{-\eta^-_{m,n}(\beta,k)}
\sin(\psi_{m,n}+\eta)d\eta
\Biggr)
\end{eqnarray}
where the coordinates $\beta$ on proj$(\what\caln^{t,1}_{m,n})$ and on
proj$(\what \caln^{t,2}_{m,n})$ (or on proj$(\caln^{t,3}_{m,n})$ and on
proj$(\what\caln^{t,4}_{m,n})$) are naturally identified.  The orbit
$(\beta,\psi_{m,n}+\eta)$ escapes not later than at the time $t$ if there
exists such integer $r$,
$0\le r\le [\frac t{2\sin(\frac mn\, \pi-\eta)}]$ that
proj~$T^r(\beta,\psi)\in H_1\cup H_2$.  

Denote $\caln^t_{\beta_0,\psi_{m,n,\ell}}=\caln^t_{\psi_{m,n,\ell}}
\cap\{(\beta,\psi):\beta=\beta_0\}$, $0\le \ell <n$.  Then, in view of
Lemma~\ref{lem4} and Theorem~\ref{thm1new}
$$\caln^t_{\beta,\psi_{m,n,\ell}}=\{(\beta,\psi):\psi_{m,n}+
\eta^-_{m,n}\le\psi\le\psi_{m,n}+\eta^+_{m,n}\},$$
where $\eta^-_{m,n}=\eta^-_{m,n}(\beta,\ell)$, $\eta^+_{m,n}=
\eta^+_{m,n}(\beta,\ell)$, $0\le \ell <n$, but we will often drop the
dependence on $\beta$ and $\ell$ to simplify notations.

Clearly $t/2\cos (\psi_{m,n}+\eta^-_{m,n})$ and 
$t/2\cos (\psi_{m,n}+\eta^+_{m,n})$ are both integers.
Moreover, for any $\eta >\eta^+_{m,n}(\eta>\eta^-_{m,n})$ the
escape occurs before the time $t$.
Therefore, either\newline 
proj~$T^{[t/2\cos(\psi_{m,n}+\eta^-_{m,n})]}
(\beta,\psi_{m,n}+\eta^-_{m,n})=
(0,\psi_{m,n}+\eta^-_{m,n})$ or\newline
proj~$T^{[t/2\cos(\psi_{m,n}+ \eta^-_{m,n})]}
(\beta,\psi_{m,n}+\eta^-_{m,n})=
(\theta,\psi_{m,n}+\eta^-_{m,n})$.\newline
Correspondingly, either \newline
proj~$T^{[t/2\cos(\psi_{m,n}+\eta^+_{m,n})]}
(\beta,\psi_{m,n}+\eta^+_{m,n})=(\udel,\psi_{m,n}+\eta^+_{m,n})$
or\newline
 proj~$T^{[t/2\cos(\psi_{m,n}+ \eta^+_{m,n})]}
(\beta,\psi_{m,n}+\eta^+_{m,n})=
(\udel+\theta,\psi_{m,n}+\eta^+_{m,n})$.

Denote by $\rho^+_{\beta,m,n}$ ($\rho^-_{\beta,m,n}$) the distances
(the lengths of the arcs) from $\beta$ to the ends (of the closest
to $\beta$) connected components of proj$(\caln^t_{\psi_{m,n}})$.
Then, according to \eref{eq1}, we have
\begin{eqnarray}\label{eq7}
\rho^+_{\beta,m,n}=(2\eta^+_{m,n}(\beta,\ell))r_\ell (\beta,m,n)\nonumber\\
\rho^-_{\beta,m,n}=(2\eta^-_{m,n}(\beta,\ell))r_\ell (\beta,m,n)
\end{eqnarray}
where $r_\ell=r_\ell(\beta,m,n)$ is a positive integer.  Analogously
\begin{eqnarray}\label{eq8}
r_\ell(\beta,m,n)2\sin\left(\frac mn\, \pi+\eta^+_{m,n}(\beta,\ell)\right)=t
\nonumber\\
r_\ell(\beta,m,n)2\sin\left(\frac mn\, \pi-\eta^-_{m,n}(\beta,\ell)\right)=t
\end{eqnarray}
Recall that $n\le \left[\frac{2\pi}\udel\right]$.  Let $t>\frac{4\pi}\udel$.
Then $r_\ell(\beta,m,n)=K(\beta,m,n)n+\ell$
where $0\le \ell\le n-1$.
It easily follows from Lemmas~\ref{lem2} and \ref{lem4} that
$K(\beta_1, m,n)=K(\beta_2,m,n)$ for any $\beta_1,\beta_2\in 
\mbox{proj}(\caln^t_{\psi_{m,n,\ell}})$, $0\le\ell <n$. Therefore,
in what follows we will write $K(m,n)$ instead of $K(\beta,m,n)$.

  From \eref{eq7} and \eref{eq8} we get after some tedious but
elementary computations
\begin{equation}\label{eq9}
1-\frac2{K(m,n)-1}<\frac{\eta^+_{m,n}(\beta,j)}{\eta^+_{m,n}(\beta,i)}<1
+\frac3{K(m,n)-1}\end{equation}
\begin{equation*}
1-\frac2{K(m,n)-1}\le\frac{\eta^-_{m,n}(\beta,j)}{\eta^-_{m,n}(\beta,i)}<1
+\frac3{K(m,n)-1}\end{equation*}
for any pair of integers $i,j$, $0\le i\le n-1$, $0\le j\le n-1$.
Another easy estimate gives that for any $\ell$, $0\le \ell <n$
\begin{equation}\label{eq10}
K(m,n)-1<\frac t{2n\sin\left(\frac mn\, \pi+\eta^+_{m,n}(\beta,j)\right)}
<K(m,n)+1,\end{equation}
\begin{equation*}
K(m,n)-1<\frac t{2n\sin\left(\frac mn\, \pi-\eta^-_{m,n}(\beta,j)\right)}
<K(m,n)+1.\end{equation*}
Finally \eref{eq2}, \eref{eq6}, \eref{eq7}, \eref{eq8}, \eref{eq9}, and
\eref{eq10} imply the following estimate of $\mu (\caln_t)$
\begin{eqnarray}\label{eq11}
\fl \frac1{4\pi}\sum_{\stackrel{\scriptstyle m<n,(m,n)=1,}
{n<\left[\frac{2\pi}\udel\right]}}
\left(1-\frac2{K(m,n)-1}\right)
\frac{n\left[g\left(\frac{2\pi}n-\theta'-\ep\right)+g(\theta'-\ep)\right]}t
\, \sin^2\frac mn\, \pi\nonumber\\
\fl<\mu(\caln_t)
<\frac1{4\pi}\sum_{\stackrel{\scriptstyle m<n,(m,n)=1,}
{n<\left[\frac{2\pi}\udel\right]}}
\left(1+\frac3{K(m,n)-1}\right)
\frac{n\left[g\left(\frac{2\pi}n-\theta'-\ep\right)+g(\theta'-\ep)\right]}t
\sin^2\frac mn\, \pi\nonumber\\
\end{eqnarray}
where $g(x)=\cases{x^2, &{$x>0$}\\0, & {otherwise}.}$

\begin{thm}\label{thm1}
\begin{eqnarray}\label{eq12}
P_\infty(\theta,\udel)&=\lim_{t\to\infty} t\mu (\caln_t)=
\frac1{8\pi}\sum^{[2\pi/\udel]}_{n=1}n(\phi(n)-\mu (n))\nonumber\\
&\qquad 
\times \left[ g\left(\frac{2\pi}n-\theta'-\udel\right)+g(\theta'-\udel)\right],
\end{eqnarray}
where $\phi (n)$ is the Euler  function and $\mu (n)$ is the 
M\"obius function.
\end{thm}

\pf Clearly
$\lim_{t\to\infty}
K(m,n)=\infty$.  We apply now the Ramanujan identity \cite{5}
\begin{equation}\label{eq13}
\sum^{n-1}_{\stackrel{\scriptstyle m= 0}{(m,n)=1}} 
e^{2\pi im/n}=\mu (n)\end{equation}
Then
\begin{eqnarray*}
\sum^{n-1}_{\stackrel{\scriptstyle m=0}{(m,n)=1}} 
\sin^2\left(\frac{\pi m} n\right)&=-\frac14
\sum^{n-1}_{\stackrel{\scriptstyle m=0}{(m,n)=1}} 
(e^{\pi im/n}-e^{-\pi im/n})^2=\\
&=-\frac14\left(2\mu (n)-2\sum^{n-1}_{\stackrel{\scriptstyle m=0}{(m,n)=1}} 
1\right)=
\frac{\phi (n)-\mu (n)}2,\end{eqnarray*}
which together with \eref{eq11} imply \eref{eq12}.

\section{The limit of small holes}

The function $P_\infty (\theta,\udel)$ is piecewise smooth with respect to
each of $\theta$ and $\udel$.  The sum in \eref{eq12} is finite, which
becomes infinite when $\udel \to 0$. 

We will study the limiting behavior of $P_\infty(\theta,\udel)$ as
$\udel\to 0$.  First, we change $[2\pi/\udel]$ onto $\infty$ in the upper
limit of $\sum$ in \eref{eq12} just by formally adding to the sum terms
identically equal zero for finite $\udel$. 

Consider now the Mellin transform
\begin{equation}\label{eq14}
\tilde P_\theta (s)=\int^\infty_0 P_\infty(\theta,\udel)\udel^{s-1}
d\udel\end{equation}
The transform $\tilde P_\theta (s)$ exists if the integral 
$\int^\infty_0 |P_\infty(\theta,\udel)|\udel^{k-1}d\udel$ is bounded
for some $k>0$.  It is certainly the case because $P_\infty(\theta,\udel)=0$
when $\udel >\pi$, and therefore \eref{eq14} converges for sufficiently large 
$s$.  

Then the inverse transform
\begin{equation}\label{eq15}
P_\infty (\theta,\udel)=\frac1{2\pi i}\int^{c+i\infty}_{c-i\infty}
\udel^{-s} \tilde P_\theta (s)ds\end{equation}
also exists if $c>k$, i.e., if $c$ is greater than the real parts of all the 
poles of $\tilde P_\theta(s)$.  

Write $\theta'=\frac{2\pi}n\left\{\frac{n\theta}{2\pi}\right\}$, where
$\{x\}$ is the fractional part of $x$.  Then
\begin{eqnarray*}
&\fl\tilde P_\theta (s)=\frac1{8\pi}\sum^\infty_{n=1}n(\phi (n)-\mu(n))
\Biggl[ \int^\infty_0\udel^{s-1}\left(\frac{2\pi}n-\theta'-\udel\right)^2\\
&\quad + \int^\infty_0 \udel^{s-1}(\theta'-\udel)^2d\udel\Biggr]\\
&= \frac1{4\pi}\sum^\infty_{n=1}n(\phi (n)-\mu (n))
\left(\frac{\left(\frac{2\pi}n-\theta'\right)^{s+2}+\theta'^{s+2}}
{s(s+1)(s+2)}\right)\\
&= \frac12\frac{(2\pi)^{s+1}}{s(s+1)(s+2)}\sum^\infty_{n=1}
\frac{\phi (n)-\mu(n)}{n^{s+1}}
\left[\left(1-\left\{\frac{n\theta}{2\pi}\right\}\right)^{s+2}+
\left\{\frac{n\theta}{2\pi}\right\}^{s+2}\right]
\end{eqnarray*}

If $s$ is sufficiently large then the convergence of the series
\eref{eq14} is uniform in $\udel$.  Therefore we can interchange
the sum and the integral in \eref{eq14}, and then integrate over 
$\udel$.
The result is
\begin{eqnarray}\label{eq16}
P_\infty(\theta,\udel)&=\frac1{2\pi i}\int^{c+i\infty}_{c-i\infty}
\frac{ds\udel^{-s}(2\pi)^{s+1}}{2s(s+1)(s+2)}
\sum^\infty_{n=1}\frac{\phi (n)-\mu (n)}{n^{s+1}}\nonumber\\
&\qquad \times\left[\left(1-\left\{\frac{n\theta}{2\pi}\right\}\right)^{s+2}
+\left\{\frac{n\theta}{2\pi}\right\}^{s+2}\right].
\end{eqnarray}

\section{Rational angles between holes}
In what follows we suppose that the angle $\theta$ between the holes
$H_1$ and $H_2$ is a rational multiple of $\pi$, i.e.,
$\theta=2\pi\, \frac rq$, $(r,q)=1$.  In particular, one gets a single
hole case by letting $r=0$, $q=1$.

For any positive integer $a$ consider the sum
\begin{equation}\label{eq17}
\sum_{n\equiv a(\bmod\; q)}\frac{\phi (n)-\mu (n)}{n^{s+1}}
\end{equation}
Transform now \eref{eq17} by dividing all the terms through by the
greatest common divisor $b=(a,q)$.  Then we get
\begin{equation}\label{eq18}
\sum_{n'\equiv a'(\bmod\; q')} \frac{\phi (bn')-\mu (bn')}{(bn')^{s+1}}
\end{equation}
where $n'=n/b$, $a'=a/b$, $q'=q/b$, where $(a',q')=1$. 

To make the paper self-contained we recall now some facts about the
Dirichlet characters (see, e.g., \cite{5} for more details).
Dirichlet's characters to the modulus $q$ are multiplicative functions
$\chi (n)$ of an integer variable $n$ which are periodic with period $q$.
The conjugacy classes modulo $q$ which are coprime to $q$ form an
abelian group under multiplication.

It is easy to see that the order of this group equals $\phi (q)$.
Besides it is a finite abelian group. Therefore it has $\phi (q)$ irreducible
representations $\chi (n)$ where $(n,q)=1$.  The characters
$\chi (n)$ are in this case the complex roots of unity, i.e.,
$\chi (m)\chi (n)=\chi (mn)$.  This definition is extended by setting
$\chi (n)=0$, if $(n,q)>1$.

By the orthogonality relation \cite{5}
\begin{equation}\label{eq19}
\frac1{\phi (q)}\sum_\chi\bar\chi (a)\chi (n)=\delta_{a,n}\end{equation}
where $\delta_{a,n}=1$, if $a\equiv n(\bmod\; q)$, zero otherwise,
and $\bar x$ denotes a complex conjugate to a number $x$.

By inserting \eref{eq19} into \eref{eq18} we get
\begin{equation}\label{eq20}
\sum_{n\equiv a(\bmod q)}\frac{\phi (n)-\mu (n)}{n^{s+1}}
=\frac1{\phi (q')}\sum_\chi \bar \chi (a')\sum^\infty_{n'=1}
\chi (n')\frac{\phi (bn')-\mu (bn')}{(bn')^{s+1}}\end{equation}
Let $n'=\prod_p p^{\alpha_p} $ be the decomposition of $n'$ into
prime factors.  Then $\chi (n')=\prod_p\chi (p)^{\alpha_p}$.
Furthermore
\begin{equation}\label{eq21}
\mu (bn')=\cases{\mu (b)\prod_p(-1)^{\alpha_p} & if $bn'$ is square free\\
0 & otherwise,}\end{equation}
\begin{equation}\label{eq22}
\phi(bn')=\phi (b)\prod_{p\mid n', p\mid b}(1-p^{-1}),
\end{equation}
where $\alpha_p=0$ if $p\mid b$.
Farther
\begin{equation}\label{eq23}
\sum^\infty_{n=1}\frac{\mu (n)}{n^{s+1}}=\prod_p (1-p^{-s-1})=
(\zeta (s+1))^{-1},\end{equation}
where $\zeta(s)$ is the Riemann zeta function. 
Now by making use of the M\"obius transform we get
\begin{equation}\label{eq24}
\sum^\infty_{n=1} \frac{\phi (n)}{n^{s+1}}=(\zeta (s+1))^{-1}
\sum^\infty_{n=1}\frac1{n^s}=\frac{\zeta (s)}{\zeta (s+1)}\end{equation}
Therefore
\begin{equation}\label{eq25}
\sum_n\frac{(\phi (n)-\mu (n))}{n^{s+1}}=\frac{\zeta (s)-1}{\zeta (s+1)}
\end{equation}
Analogously
\begin{equation}\label{eq26}
\sum_n\frac{\chi (n)(\phi (n)-\mu (n))}{n^{s+1}} = 
\frac{L(s,\chi )-1}{L(s+1,\chi)}\end{equation}
Finally we have
\begin{eqnarray}\label{eq27}
\sum_{n\equiv a(\bmod q)}\frac{\phi (n)-\mu (n)}{n^{s+1}}\nonumber\\
\qquad=\frac1{b^{s+1}\phi (q')} \sum_\chi 
\frac{\bar \chi (a')(\phi (b)L(s,\chi )-\mu (b))}{L(s+1,\chi )\prod_{p\mid b}
(1-\chi (p)p^{-s-1})}\end{eqnarray}
where the characters are taken modulo $q'$ and 
\begin{equation}\label{eq28}
L(s,\chi)=\sum^\infty_{n=1}\frac{\chi (n)}{n^s} = \prod_p
\left(1-\frac{\chi (p)}{p^s}\right)^{-1}\end{equation}
is the Dirichlet $L$ function.

If $q'=1$ then $L(s,\chi)$ reduces to the Riemann zeta function $\zeta (s)$.
For each $q'$ there is a trivial character $\chi (a')$ that assumes the 
value 1 for all $a'$ coprime to $q'$.  Therefore
\begin{equation}\label{eq29}
L(s,1)=\zeta (s)\prod_{p\mid q'} (1-p^{-s}).\end{equation}

Let
\begin{eqnarray}\label{eq31}
\tilde P_{r/q}(s)&=\frac{(2\pi)^{s+1}}{2s(s+1)(s+2)}
\sum^q_{a=1}\frac{\left(1-\left\{\frac{ar}q\right\}\right)^{s+2}
+\left\{\frac{ar}q\right\}^{s+2}}{b^{s+1}\phi(q')}\nonumber\\
&\qquad\times\sum_\chi \frac{\bar\chi (a')(\phi (b)L(s,\chi )-\mu (b))}
{L(s+1,\chi)\prod_{p\mid b}(1-\chi (p)p^{-s-1})},\end{eqnarray}
where, as above, $b=(a,q)$, $a'=a/b$, $q'=q/b$ and the characters are
taken $\bmod\; q'$.  We note that odd characters (i.e. $\chi(-1)=-1$) and their
$L$ functions in the above expression cancel.

The function $\tilde P_{r/q}(s)$ has poles at $s=0$,
$s=-1$, $s=-2$, at zeros of $L(s+1,\chi)$ and at poles of $L(s,\chi)$.
Dirichlet's function $L(s+1,\chi)$ with even $\chi$ has trivial zeros
at $s=-(2m+1)$, where $m=1,2,\dots$ \cite{6}.  All other (nontrivial) zeros 
of $L(s+1,\chi)$ have real part $Re\,s=-1/2$ assuming that the 
extended Riemann hypothesis that is concerned with the Dirichlet
functions \cite{5} is correct.

\section{The simplest placements of two holes and the Riemann hypothesis}
In this section we consider several specific values of $q$, 
when a number of characters does not exceed 2, i.e., $\phi (q)\le 2$.
There are thus five such values $q=1,2,3,4$ and 6. In all these cases
the only even character is the trivial character, so the function
$\tilde P_{r/q}(s)$ in \eref{eq31} contains the Riemann zeta function and
no other $L$ functions.

Below there are two tables.  The first
table gives the exact expressions for the function
$\tilde P_{1/q} (s)$ for $q=1,2,3,4,6$ and $r=1$.  The second
table contains the corresponding residues.
(Except for the last lines for $q=6$ these two tables were 
published in \cite{1}.)

We now list some properties of the Riemann zeta function that will
be needed in what follows.
Riemann proved that $\zeta$-function satisfies to the following
functional equation
\begin{equation}\label{eq32}
\ugam\left(\frac s2\right)\pi^{-s/2}\zeta (s)=
\ugam\left(\frac{1-s}2\right)\pi^{-(1-s)/2}\zeta (1-s)
\end{equation}

\begin{table}
\begin{indented}
\item[]\begin{tabular}{|c|c|}\hline
$q$&$\tilde{P}(s)\phantom{\Big)}$\\\hline
1&$\displaystyle\frac{(2\pi)^{s+1}(\zeta(s)-1)}{2s(s+1)(s+2)\zeta(s+1)}$\\[1em]
2&$\displaystyle\frac{\pi^{s+1}\zeta(s)}{s(s+1)(s+2)\zeta(s+1)}$\\[1em]
3&$\displaystyle\frac{(2\pi/3)^{s+1}(3^s(7\zeta(s)+2^{s+2}(\zeta(s)-1)+2)
-\zeta(s)(2^{s+2}+1)}{2s(s+1)(s+2)(3^{s+1}-1)\zeta(s+1)}$\\[1em]
4&$\displaystyle\frac{(\pi/2)^{s+1}(2^s(13\zeta(s)+3^{s+2}(\zeta(s)-1)+3)
-\zeta(s)(3^{s+2}+5))}{4s(s+1)(s+2)(2^{s+1}-1)\zeta(s+1)}$
\\
&\raisebox{-8pt}[-8pt][8pt]{$\scriptstyle[(\pi/3)^{s+1}(6^s+8.12^s-25.30^s+
(1-3.2^s-13.3^s-8.4^s$}\\6&$\scriptstyle
+25.5^s+27.6^s-25.10^s+8.12^s-25.15^s+25.30^s)\zeta(s))
]$\\&\raisebox{4pt}[4pt][-4pt]{$\scriptstyle\times[
2s(s+1)(s+2)(2^{s+1}-1)(3^{s+1}-1)\zeta(s+1)]^{-1}
$}\\\hline
\end{tabular}
\caption{The function $\tilde{P}_{r/q}(s)$ (Eq.~\eref{eq31})
for $q=1,2,3,4,6$ and $r=1$.
\label{t:P(s)}}
\end{indented}
\end{table}
\vskip-.5\baselineskip
\begin{table}
\begin{indented}
\item[]\begin{tabular}{|c|cccc|}\hline
$q$&\multicolumn{4}{c|}{$s$}\\
&1&$-1$&$-2$&$-3$\\\hline
1&2&$-\frac{13}{12}$&$\frac{3}{2\pi}$&$\frac{119}{5760\pi^2\zeta'(-2)}$\\
2&1&$-\frac{1}{6}$&0&$-\frac{1}{720\pi^2\zeta'(-2)}$\\
3&1&$-\frac{1}{4}-\frac{5\ln2}{9\ln3}$&$\frac{3}{4\pi}$&
$\frac{49}{5120\pi^2\zeta'(-2)}$\\
4&1&$-\frac{1}{3}-\frac{11\ln3}{16\ln 2}$&$\frac{3}{\pi}$&
$\frac{109}{1620\pi^2\zeta'(-2)}$\\
&&\raisebox{-8pt}[-8pt][8pt]{$
\scriptstyle[5\ln5(10\ln3-7\ln5)+\ln2(55\ln5-76\ln3)$}&&\\
6&1&$\scriptstyle+(10\ln5-8\ln2)
(7\ln\udel+12\zeta'(-1))]$&$-\frac{3}{2\pi}$&
$-\frac{79}{6400\pi^2\zeta'(-2)}$\\
&&\raisebox{4pt}[4pt][-4pt]{$\scriptstyle\times[72\ln2\ln3]^{-1}$}&&
\\\hline
\end{tabular}
\caption{Some residues of $\tilde{P}_{r/q}(s)\udel^{-s}$
given in table~\protect\ref{t:P(s)} divided by the factor $\udel^{-s}$.
The $\ln\udel$ appears for $q=6$ due to a double pole at $s=-1$.
There are also poles for further negative odd $s$, and along the
critical line ${\cal R}e\;s=-1/2$.
\label{t:poles}.}
\end{indented}
\end{table}

\noindent
which can also be written in the following (nonsymmetric) form
\begin{equation}\label{eq33}
\zeta (1-s)=2^{1-s}\pi^{-s}\cos(\frac\pi 2\, s)\ugam (s)\zeta (s)
\end{equation}
where $\ugam (s)$ is the gamma function.

It is well known also that $\zeta (0)=-\frac12$, $\zeta (-2m)=0$,
\begin{equation}\label{eq34}
\zeta (1-2m)=\frac{(-1)^mB_{2m}}{2m}
\end{equation}
where $m=1,2\dots$, and $B_1,B_2,\dots$ are Bernoulli numbers.
 The following approximating formula holds
for the Bernoulli numbers
\begin{equation}\label{eq35}
B_{2m}\sim (-1)^{m-1}4\sqrt{\pi m}\left(\frac m{\pi e}\right)^{2m}.
\end{equation}

Another well known fact is that $\ugam (s)$ has poles of order 1
at $s=-m$ for all integers $m>0$ and
\begin{equation}\label{eq36}
Res (\ugam,-m)=\frac{(-1)^m}{m!}\end{equation}
where $Res (f,a)$ denotes the residue of the function $f(x)$ at
the point $x=a$.

\begin{lem}\label{lem5}
Let $q=1,2,3,4$ or $6$, then
$$\sum_j\mathop{Res}\limits_{s=s_j}(\tilde P_{1/q}(s)\udel^{-s})<C\udel
|\ln \udel| ,$$
where $C>0$ is a constant and the sum is taken over $s_j=-1,-2$ and over
all trivial zeros of $\zeta (s+1)$, i.e., over all 
odd negative integers $m\le -3$.
\end{lem}

By combining \eref{eq33}, \eref{eq34}, \eref{eq35}, \eref{eq36}
with \eref{eq31} and Stirling formula we obtain that
$$\Bigl| \sum_j \mathop{Res}\limits_{s=s_j}(\tilde P_{1/q}(s))\Bigr|
<\sum_j\bigl| \mathop{Res}\limits_{s=s_j}(\tilde P_{1/q}(s))\bigr|<
\mbox{const}$$
and Lemma~\ref{lem5} follows.

\begin{rem}\rm
The (extra) factor $|\ln \udel|$ appears in the right hand side of the 
above estimate because of the double pole at $s=-1$ for $q=6$.
For $q=1,2,3,4$ this factor is not needed.
\end{rem}

It is well known \cite{6} that for $\tau\ge \tau_0 >0$ uniformly in 
$\sigma$ the following estimates hold
\begin{equation}\label{eq36new}
\zeta (\sigma+i\tau)=O\cases{
1, & $\sigma\ge 2$\\
\log\tau, & $1\le \sigma\le 2$\\
\tau^{(1-\sigma)/2}\log\tau, & $0\le \sigma\le 1$\\
\tau^{1/2-\sigma}\log\tau, & $\sigma\le 0.$}
\end{equation}

We will consider now nontrivial zeros of the Riemann zeta function
located at the critical strip $0<\sigma <1$.
Let $N(t)$ denotes the number of zeros of $\zeta (s)=\zeta (\sigma + it)$
in the region $\{(\sigma, t):0<\sigma <1$, $0<t\le T\}$ of the critical
strip.  Then \cite{6}
\begin{equation}\label{eq37}
N(T)=\frac T{2\pi}\log\frac T{2\pi} - \frac T{2\pi} + \frac 78
+O(\log T)\end{equation}
We now assume that the Riemann hypothesis (RH) is correct and use
several of its well known \cite{6} consequences.

Let $S(t)$ denote the multiplicity of the complex zero $S=\frac12 + it$ of
$\zeta (s)$.  Then on RH
\begin{equation}\label{eq40}
S(t)=O\left(\frac{\log t}{\log\log t}\right)\end{equation}
We will construct an infinite sequence of contours $C_n$ over which the
integration of the function $\tilde P_{r/q}(s)$ (see \eref{eq31}) will be
performed in what follows.  Each contour in this sequence will contain two
vertical segments
$$I_n(k_0)=\{s=\sigma+i\tau:\sigma=k_0, -a_n\le\tau\le a_n\}$$
and $I'_n = \{s=\sigma+i\tau:\sigma=-b_n, -a_n\le\tau\le a_n\}$, where
$b_n=2n$, and two horizontal segments
$I^+_n=\{s=\sigma+i\tau:-b_n\le\sigma\le k_0,\tau=a_n\}$ and
$I^-_n=\{s=\sigma+i\tau:-b_n\le\sigma\le k_0,\tau=-a_n\}$.
We choose $k_0$ large enough to ensure uniform convergence of the 
Mellin transform \eref{eq14} and assume from now on that 
$\udel<1$..  

\begin{lem}\label{lem6new}
There exists an infinite sequence of contours $C_n$ with $a_n\to\infty$ as
$n\to\infty$ such that
$$\lim_{n\to\infty}\int_{I'_n\cup I^+_n\cup I^-_n}\tilde P_{r/q}(s)
\udel^{-s}ds=0$$
for any entry in the {\rm Table~1} (i.e., for $q=1,2,3,4,6$).  
\end{lem}

\pf One gets from the Table~1 that $s=1$ is the only pole
of $\tilde P_{r/q}(s)$.  Therefore we will assume that
$k_0>1$. Besides, it is easy to see from the Table~1 that it is enough to 
consider instead of $\tilde P(s)$ the function
$$\hat P(s)=\frac{(2\pi)^{s+1}}{2q^{s+1}s(s+1)(s+2)}
\frac{\zeta (s)}{\zeta (s+1)}\,.$$
We will start with the vertical segment $I'_n$.  By making use of the
functional equation \eref{eq33} we get
\begin{equation}\label{eq39new}
\frac{\zeta (s)}{\zeta (s+1)} = -\frac s2
\frac{\sin\frac\pi2 s}{\sin \frac{(\pi+1)s}2 }
\frac{\zeta(1-s)}{\zeta (-s)}\end{equation}

Because of the relation $\overline {\zeta (\bar s)}=\zeta (s)$, where
$\bar s$ denotes the complex conjugate to the complex number $s$, it is 
enough to consider $\int_{I^+_n}\hat P(s)\udel^{-s}ds$.  The estimates
for $\int_{I^-_n}\hat P(s)\udel^{-s}ds$ are quite analogous.

Clearly, the major problem with estimating of these integrals is
caused by the zeros of $\zeta (s+1)$ in the denominator of 
$\hat P(s)$. Therefore we partition the horizontal segment 
$I^+_n$ into the union of three segments
\begin{eqnarray*}
I^+_{n,1} &= \{s=\sigma +i\tau: -2n\le \sigma\le -1,\tau=a_n\},\\
I^+_{n,2} &= \{s=\sigma +i\tau: -1< \sigma< 0,\tau=a_n\},\quad\mbox{and}\\
I^+_{n,3} &= \{s=\sigma +i\tau: 0\le \sigma\le k_0,\tau=a_n\}.\end{eqnarray*}

Observe that the lengths of $I^+_{n,2}$ and $I^+_{n,3}$ are constants.
These sets are not exactly defined though because the values of $a_n$ are
not specified yet.  We will make a choice of $a_n$ now.  To do that we
assume the validity of the Riemann hypothesis.
Then~\cite{6} each interval $(\tau,\tau+1)$ on the critical line 
$\sigma=\frac12$ contains a value of $\tau$ such that
\begin{equation}\label{eq40new}
|\zeta (s)| >\exp\left(-A_1\frac{\log\tau\log\log\log\tau}
{\log\log\tau}\right)\end{equation}
where $A_1$ is an absolute constant.

We now choose $a_n$ so that $0<a_n\le n+1$ and 
$\zeta\left(\frac12+ia_n\right)$ satisfies \eref{eq40new}.
The next two estimates that we will use also hold on the assumption
that the Riemann hypothesis is correct~\cite{6}.

The first fact is that for any sufficiently small  $\ep >0$
\begin{equation}\label{eq41new}
\zeta (s)=O(\tau^\ep)\quad\mbox{and}\quad
\frac1{\zeta (s)} = O(\tau^\ep)\end{equation}
if $\frac12 +\frac1{\log\log\tau} \le\sigma$.
Clearly \eref{eq39new} together with \eref{eq41new} gives
that $\lim_{n\to\infty}\int_{I'_n}\hat P(s)\udel^{-s}ds=0$.

Another estimate works in the vicinity of the critical line
$\frac12\le\sigma\le \frac12+\frac1{\log\log\tau}$ \cite{6}
\begin{equation}\label{eq42new}
\log|\zeta (s)|>\frac{-A_2\log \tau}{\log\log\tau}
\log\left\{\frac2{\left(\sigma-\frac12\right)\log\log\tau}\right\},
\end{equation}
where $A_2$ is another absolute constant.

Now \eref{eq40new}, \eref{eq41new}, \eref{eq42new} and the last
(fourth) estimate in \eref{eq36new} applied for $-1\le\sigma\le 0$ imply
that $\lim_{n\to\infty} \int_{I^+_{n,2}}\hat P(s)\udel^{-s}ds=0$.

In the estimates of the integrals of $\hat P(s)$ over $I^+_{n,1}$ and
$I^+_{n,3}$ the term $s(s+1)(s+2)$ in the denominator ensures the needed
results.
Indeed, consider first $I^+_{n,3}$.  Then \eref{eq41new} implies that
$\int^+_{I_{n,3}}\hat P(s)\udel^{-s}ds\mathop{\longrightarrow}\limits_{n\to 0}
0$.  For $I^+_{n,1}$ we will again use the trick with the reflection
(the functional equation).  Then the relation \eref{eq39new} together 
with the estimates \eref{eq41new} ensures that on
$I^+_{n,1}\hat P(s)$ satisfies the inequality
$|\hat P(s)|<|n|^\gamma$, where $0<\gamma <1$.  The length of 
$I^+_{n,1}$ equals $2n$.  Thus $\lim_{n\to\infty}\int_{I^+_{n,1}}
\hat P(s)\udel^{-s}ds=0$ and Lemma~\ref{lem6new} follows.

\begin{lem}\label{lem6}
Assume that RH is correct.  Let $q=1,2,3,4$ or $6$, $(r,q)=1$, then
for any $\alpha>0$
\begin{equation}\label{eq42}
C_2\udel^{1/2}<\sum_j \mathop{Res}\limits_{s=\frac12+i\tau_j}
(\tilde P_{r/q}(s)\udel^{-s})<C_1\udel^{1/2-\alpha},\end{equation}
where the sum is taken over all nontrivial zeros of $\zeta (s)$
and $C_1,C_2>0$ are some constants.
\end{lem}

\pf Observe, at first, that all expressions for $\tilde P_{r/q}(s)$ in the
Table~1 have $\zeta (s+1)$ in denominators.
Therefore all poles outside the real line correspond to zeros of 
$\zeta (s+1)$ and under RH are located on the line 
$s=-\frac12 + i\tau$, $-\infty < \tau<\infty$.  

It is
easy to see that residue at each zero of $\zeta (s+1)$ with multiplicity
$m$ results in the extra factor $(\ln \udel)^{m-1}$ in the expression for
the residue at this zero.

By making use of \eref{eq40} we get for any fixed $\udel$ and
for sufficiently large $\tau$
\begin{equation}\label{eq43}
|\log\udel|^{S(\tau)}<\tau^\delta\end{equation}
where $0<\delta <1/2$.

  From \eref{eq37}, \eref{eq40}, \eref{eq41new} and \eref{eq43}
we have for sufficiently large $T$
\begin{eqnarray}\label{eq44}
\Bigl|\sum^\infty_{j=1}\mathop{{Res}}\limits_{s=\frac12+i\tau_j}
(\tilde P_{1/q}(s))\Bigr| &= 
\Bigl|\sum^\infty_{n=0} \sum_{n\le \tau_j\le n+1}
\mathop{Res}\limits_{s=\frac12+i\tau_j}
(\tilde P_{1/q}(s))\Bigr|\nonumber\\
& \le
\sum^\infty_{n=0}\Bigl|\sum_{n\le \tau_j<n+1}{Res}
(\tilde P_{1/q}(s))\Bigr|\nonumber\\
& \le\sum^\infty_{n=0}
\frac{(n+1)^{\del+\ep}}{n^3}<\infty
\end{eqnarray}
Observe
that convergence (or divergence) of the series in 
\eref{eq42new} is not influenced by the
constant factor $\udel^{1/2}$.
Thus Lemma~\ref{lem6} follows.  

  From Lemmas~\ref{lem5}, \ref{lem6new} and \ref{lem6} follows

\begin{thm}\label{thm3new}
Consider a billiard in the unit circle with two holes $[0,\udel]$ and
$\left[ 2\pi\frac rq\,, 2\pi \frac rq+\udel\right]$, where $q=1,2,3,4$ and
$6$, $0<r<q$ are integers, $(r,q)=1$.  If $t>f(t)\udel^{-1}$, where
$f(t)>0$ and $\lim_{t\to\infty} f(t)=\infty$, then
\begin{equation}\label{eq46new}
P_\infty\left(\frac rq,\udel\right)=\lim_{t\to\infty}
\mu (\caln_t)=\sum_k\mathop{Res}\limits_{s=s_k}
\tilde P_{r/q}(s)\udel^{-s}\end{equation}
where summation is taken over all residues of the function
$\tilde P_{r/q}(s)$.
\end{thm}

\pf Consider the sequence of contours $C_n$ constructed in the proof of
Lemma~\ref{lem6new}.  Then 
$$\lim_{n\to\infty}\int_{C_n}\tilde P_{r/q}(s)\udel ^{-s}ds=
\int^{k_0+i\infty}_{k_0-i\infty}\tilde P_{r/q}(s)\udel^{-s}ds.
$$
Farther applying the residue theorem and Lemmas~\ref{lem5} and \ref{lem6}
to the integral $\int_{C_n}\tilde P_{r/q}(s)\udel^{-s}ds$ and letting
$n\to\infty$ we obtain that $\int^{k_0+i\infty}_{k_0-i\infty}
\tilde P_{r/s}(s)\udel^{-s}ds=\sum_k\mathop{Res}\limits_{s=s_k}
\tilde P_{r/q}(s)\udel^{-s}<\infty$.  The relation
\eref{eq46new} follows from \eref{eq15} and \eref{eq31}.

Combining now Theorems~\ref{thm1} and \ref{thm3new} with Lemmas~\ref{lem5},
\ref{lem6new} and \ref{lem6} and with the first row in the Table~2 we obtain

\begin{thm}\label{thm3}
Consider an open circular billiard with one hole (i.e., with two holes
of the same length $\udel$ placed on top of each other).  Let $P_1(t,\udel)$
denotes the probability that a particle will not escape till time $t$.
Then, assuming that RH is correct, for any $\alpha >0$
\begin{equation}\label{eq45}
\lim_{\udel \to 0}\lim_{t\to\infty}\udel^{\alpha-1/2}
[tP_1(t,\udel)-2/\udel]=0\end{equation}
\end{thm}

The inverse statement is also true.

\noindent{\bf Theorem~\ref{thm3}$'$}~~{\it Consider the same
open circular billiard as in Theorem~{\rm\ref{thm3}}.  Then the validity of the
relation {\rm\eref{eq45}} implies RH.}

\pf Suppose that RH is not true.  Then the Riemann zeta function
$\zeta (s)$ has at least one zero $s_0$ in the critical strip
which is outside the critical line, i.e., 
$s_0=\left(\frac12+\gamma+i\tau\right)$,
where $\gamma\ne 0$, $|\gamma|<\frac12$.  The functional equation
\eref{eq32} ensures that $\zeta (s_1)=0$, where $s_1=\frac12-
\gamma+it$.  But either $Re\; s_0<\frac12$ or $Re\; s_1<\frac12$,
which imply that instead of $\udel^{\alpha-1/2}$ in \eref{eq45}
one must have a $\udel^{\alpha-\min (Re\; s_0, Re\; s_1)}$  
unless the sum of residues of $\tilde P_{1/q}(s)\udel^{-s}$ over all zeros
in the critical strip with the real parts less than $\frac12$ is
identically (for all $0<\udel\le \udel_0)$ equal to zero.  It is easy to
see that it cannot occur.  Thus we come to the contradiction which 
proves Theorem~\ref{thm3}$'$.

Consider now circular billiard with two opposite (symmetric with respect
to the center) holes with lengths $\udel$.  Denote by $P_2(t,\udel)$
probability that the billiard particle will not escape from this circle 
till time $t$.

  From the first two rows in the Tables~1 and 2 we have that the
relation \eref{eq45} is equivalent to the statement that
\begin{equation}\label{eq46}
\lim_{\udel\to 0}\lim_{t\to\infty} \udel^{\alpha-1/2}
[tP_1(t,\udel)-2tP_2(t,\udel)]=0\end{equation}
for any $\alpha >0$.
Thus one can formulate the analogs of the Theorems~\ref{thm3} and 
\ref{thm3}$'$ by substituting \eref{eq46} instead of \eref{eq45}.
Therefore RH is equivalent to \eref{eq46} which relates asymptotics
of probabilities to escape in a circular billiard with one and
with two symmetric holes.

Certainly one can use another rows in Tables~1 and 2 to formulate
statements equivalent to RH in open billiards with two holes
places under the angles $\frac{2\pi}3$, $\frac\pi 2$ and
$\frac\pi 3$.  

Moreover, in fact the analogous statements hold for open circular
billiards with any number of $q$ holes with lengths $\udel$ which
are equally spaced over
the circle on the 
angle $\frac{2\pi}q$.

\begin{thm}\label{thm5new}
Consider an open circular billiard with $q\geq 2$ holes of the same length
$\udel$ with the centers placed at the vertices of a right convex
$q$-angle.  Let $P_q(t,\udel)$ denotes the probability that the particle
will not escape through this system of $q$ (different) holes and
$P^{(q)}_1 (t,\udel)$ denotes the probability that the particle will not
escape till time $t$ in case when all these holes are placed on top of 
each other.  Then, assuming that {\rm RH} is correct, for any 
$\alpha>0$
\begin{equation}\label{eq49new}
\lim_{\udel\to 0} \lim_{t\to \infty}\udel^{\alpha-\frac12}
t[P^{(q)}_1(t,\udel)- qP_q(t,\udel)]=0\end{equation}
\end{thm}

The proof of Theorem~\ref{thm5new} is completely analogous to the proof
of Theorem~\ref{thm3}.  Theorem~\ref{thm1new} remains valid for this system.  The formula in Theorem~\ref{thm1} becomes
\[ P_q(\udel)=\lim_{t\to\infty} t\mu(\caln_t(q))=
\frac{1}{8\pi}\sum_{n\geq 1}n(\phi(n)-\mu(n))\tilde{q}
g\left(\frac{2\pi}{n\tilde{q}}-\udel\right)
\]
where $\tilde{q}=q/gcf(n,q)$ and the sum is over $n$ for which the argument of the $g$ function is positive.  The sum is then written in terms of $\tilde{n}=n/gcf(n,q)$ with the result
\[ P_q(\udel)=\frac{1}{8\pi}\sum_{\tilde{n}=1}^{[2\pi/\udel]}\tilde{n}\phi(\tilde{n})
q^2 g\left(\frac{2\pi}{\tilde{n}}-\udel\right) \]
The Mellin transform of this is then
\begin{equation}\label{eq50new}
\tilde{P}_q(s)=\frac{(2\pi)^s\zeta (s)}{q^ss(s+1)(s+2)\zeta (s+1)}
\end{equation}
for which the case $q=2$ has already been given in Table~1.

The proof of the next statement is completely analogous to the one of 
Theorem~\ref{thm3}$'$.  

\noindent{\bf Theorem~\ref{thm5new}$'$}~~{\it Consider the same
open circular billiard with $q$ uniformly placed holes as in 
{\rm Theorem~\ref{thm3}}.  Then the relation {\rm\eref{eq49new}} implies
{\rm RH}.}

\section{Concluding Remarks}

It is quite likely that the results of this paper can be readily
generalized for Dirichlet $L$ functions.  Indeed, the main
formula \eref{eq31} for two holes escape with arbitrary
rational angles between holes explicitly involves all
$L(s,\chi)$ with even nontrivial characters $\chi$.  We conjecture
that the corresponding $\udel^{1/2}$ asymptotics for the two holes
escape is equivalent to the extended RH for even characters.  We reserve the
term generalized RH for more general $L$-functions over number fields,
elliptic curves, etc.  It seems very interesting though whether the
generalized RH is equivalent to a particular asymptotics of the escape in
some specific classes of open dynamical systems.  The most natural 
candidates for these systems are geodesic and contact flows on 
manifolds.

Another natural further problem is to compute the second order
asymptotics of the escape from the open circular billiard through two
holes placed under irrational $(\bmod\; \pi)$ angles.  The leading
order (in $\udel$) behavior in this case remains the same though 
\cite{1}.

\ack
The work of LAB was partially supported by NSF grant \#DMS-0140165
and by the Humboldt Foundation.

\end{document}